\documentclass[10pt]{article}
\begin{document}
\title{ {\bf A Characterization On Potentially $K_{2,5}$-graphic
Sequences}
\thanks{   Project Supported by  NSF of Fujian(2008J0209),
 Fujian Provincial Training
Foundation for "Bai-Quan-Wan Talents Engineering" , Project of
Fujian Education Department and Project of Zhangzhou Teachers
College.}}
\author{{ Lili Hu , Chunhui Lai}\\
{\small Department of Mathematics, Zhangzhou Teachers College,}
\\{\small Zhangzhou, Fujian 363000,
 P. R. of CHINA.}\\{\small  jackey2591924@163.com ( Lili Hu, Corresponding author)}
 \\{\small   zjlaichu@public.zzptt.fj.cn(Chunhui
 Lai )}
}

\date{}
\maketitle
\begin{center}
\begin{minipage}{4.1in}
\vskip 0.1in
\begin{center}{\bf Abstract}\end{center}
 { For given a graph $H$, a graphic sequence $\pi=(d_1,d_2,\cdots,d_n)$ is said to be potentially
 $H$-graphic if there exists a realization of $\pi$ containing $H$ as a subgraph.
 Let $K_m-H$ be the graph
obtained from $K_m$ by removing the edges set $E(H)$ where $H$ is a
subgraph of $K_m$. In this paper, we characterize  potentially
$K_{2,5}$-graphic sequences. This characterization implies a special
case of a theorem due to  Yin et al. [26].}
\par
\end{minipage}
\end{center}
 \par
 \section{Introduction}
\par
\baselineskip 14pt
    We consider finite simple graphs. Any undefined notation follows
that of Bondy and Murty $[1]$. The set of all non-increasing
nonnegative integer sequence $\pi=(d_1,d_2,\cdots,d_n)$ is denoted
by $NS_n$. A sequence $\pi\in NS_n$ is said to be graphic if it is
the degree sequence of a simple graph $G$ of order $n$; such a graph
$G$ is referred as a realization of $\pi$. The set of all graphic
sequences in $NS_n$ is denoted by $GS_n$.  Let $C_k$ and $P_k$
denote a cycle on $k$ vertices and a path on $k+1$ vertices,
respectively. Let $\sigma(\pi)$ be the sum of all the terms of
$\pi$, and let [x] be the largest integer less than or equal to $x$.
A graphic sequence $\pi$ is said to be potentially $H$-graphic if it
has a realization $G$ containing $H$ as a subgraph. Let $G-H$ denote
the graph obtained from $G$ by removing the edges set $E(H)$ where
$H$ is a subgraph of $G$.  In the degree sequence, $r^t$ means $r$
repeats $t$ times, that is, in the realization of the sequence there
are $t$ vertices of degree $r$.
\par

  Gould et al.[8] considered an extremal problem on potentially $H$-graphic sequences as follows:
  determine the smallest even
integer $\sigma(H,n)$ such that every n-term positive graphic
sequence $\pi$ with $\sigma(\pi)\geq \sigma(H,n)$ has a realization
$G$ containing $H$ as a subgraph. For $r\times s\times t$ complete
3-partite graph $K_{r,s,t}$, Yin [23] and Lai [17] independently
determined $\sigma(K_{1,1,3},n)$. Chen [6] determined
$\sigma(K_{1,1,t},n)$ for $t\geq3$, $n\geq2[{{(t+5)^2}\over4}]+3$.
For the graph $K_m-G$, Lai [16] determined $\sigma(K_4-e,n)$ and Yin
et al.[27] determined $\sigma(K_5-e,n)$. In [29], Yin determined
$\sigma(K_{r+1}-K_3,n)$ for $r\geq3$, $n\geq3r+5$. \par
  A harder question is to characterize potentially
 $H$-graphic sequences without zero terms. Yin and Li [24] gave two sufficient conditions
 for $\pi\in GS_n$ to be potentially $K_r-e$-graphic.  Luo [20] characterized potentially
 $C_k$-graphic sequences for each $k=3,4,5$. Chen [2] characterized potentially
 $C_6$-graphic sequences.  Chen et al.[3] characterized potentially
 $_kC_l$-graphic sequences for each $3\leq k\leq5$, $l=6$. Recently, Luo and Warner [21] characterized potentially
 $K_4$-graphic sequences.  Eschen and Niu [7] characterized potentially
 $K_4-e$-graphic sequences. Yin et al.[25] characterized  potentially
 $_3C_4$, $_3C_5$ and $_4C_5$-graphic sequences. Yin and Chen [28] characterized
 potentially $K_{r,s}$-graphic sequences for $r=2,s=3$ and
 $r=2,s=4$. Yin et al.[30] characterized
 potentially $K_5-e$ and $K_6$-graphic sequences. In [31], they characterized potentially
$K_6-K_3$-graphic sequences. Moreover, Yin et al.[32] characterized
potentially $K_{1,1,s}$-graphic sequences for $s=4$ and $s=5$.  Chen
and Li [5] characterized  potentially $K_{1,t}+e$-graphic sequences.
Chen [4] characterized  potentially $K_6-3K_2$-graphic sequences.
Hu, Lai and Wang [11] characterized potentially $K_5-P_4$ and
$K_5-Y_4$ -graphic sequences, where $Y_4$ is a tree on 5 vertices
and 3 leaves. Hu and Lai [9,12] characterized potentially $K_5-C_4$
and $K_5-E_3$-graphic sequences, where $E_3$ denotes graphs with 5
vertices and 3 edges. Besides, in [13,14], they characterized
potentially $K_{3,3}$, $K_6-C_6$ and $K_6-C_4$-graphic sequences.
Recently, Liu and Lai[19] characterized potentially
$K_{1,1,2,2}$-graphic sequences. Xu and Lai [22] characterized
potentially $K_6-C_5$-graphic sequences.
\par
In this paper, we characterize potentially $K_{2,5}$-graphic
sequences. This characterization implies a
 theorem due to  Yin et al. [26].
\par
\section{Preparations}\par
   Let $\pi=(d_1,\cdots,d_n)\in NS_n,1\leq k\leq n$. Let
    $$ \pi_k^{\prime\prime}=\left\{
    \begin{array}{ll}(d_1-1,\cdots,d_{k-1}-1,d_{k+1}-1,
    \cdots,d_{d_k+1}-1,d_{d_k+2},\cdots,d_n), \\ \mbox{ if $d_k\geq k,$}\\
    (d_1-1,\cdots,d_{d_k}-1,d_{d_k+1},\cdots,d_{k-1},d_{k+1},\cdots,d_n),
     \\ \mbox{if $d_k < k.$} \end{array} \right. $$
  Denote
  $\pi_k^\prime=(d_1^\prime,d_2^\prime,\cdots,d_{n-1}^\prime)$, where
  $d_1^\prime\geq d_2^\prime\geq\cdots\geq d_{n-1}^\prime$ is a
  rearrangement of the $n-1$ terms of $\pi_k^{\prime\prime}$. Then
  $\pi_k^{\prime}$ is called the residual sequence obtained by
  laying off $d_k$ from $\pi$. For simplicity, we denote $\pi_n^\prime$ by $\pi^\prime$ in this paper.
  \par
   For a nonincreasing positive integer sequence $\pi=(d_1,d_2,\cdots,d_n)$, we write $m(\pi)$ and $h(\pi)$ to denote the largest
positive terms of $\pi$ and the smallest positive terms of $\pi$,
respectively. We need the following results.
\par
    {\bf Theorem 2.1 [8]} If $\pi=(d_1,d_2,\cdots,d_n)$ is a graphic
 sequence with a realization $G$ containing $H$ as a subgraph,
 then there exists a realization $G^\prime$ of $\pi$ containing $H$ as a
 subgraph so that the vertices of $H$ have the largest degrees of
 $\pi$.
 \par
    {\bf Theorem 2.2 [18]} If $\pi=(d_1,d_2,\cdots,d_n)$ is a
 sequence of nonnegative integers with $1\leq m(\pi)\leq2$,
 $h(\pi)=1$ and even $\sigma(\pi)$, then $\pi$ is graphic. \par
{\bf Theorem 2.3 [28]} Let $n\geq6$ and $\pi=(d_1,d_2,\cdots,d_n)\in
GS_n$. Then $\pi$ is potentially $K_{2,4}$-graphic if and only if
$\pi$ satisfies the following conditions:\par (1) $d_2\geq4$ and
$d_6\geq2$;\par (2) If $d_1=n-1$ and $d_2=4$, then $d_3=4$ and
$d_6\geq3$;\par (3) $\pi\neq(4^3,2^4)$, $(4^2,2^5)$, $(4^2,2^6)$,
$(5^2,4,2^4)$, $(5^3,3,2^3)$, $(6,5^2,2^5)$, $(5^3,2^4,1)$,
$(6^3,2^6)$, $(n-1,4^2,3^4,1^{n-7})$,  $(n-1,4^2,3^5,1^{n-8})$,
$(n-2,4^2,2^3,1^{n-6})$,  $(n-2,4^3,2^2,1^{n-6})$.
\par
    {\bf Lemma 2.4 [10]} Let $n\geq4$ and $\pi=(d_1,d_2,\cdots,d_n)$ be a
 nonincreasing sequence of positive integers with even $\sigma(\pi)$. If
$d_1\leq3$ and $\pi\neq(3^3,1),(3^2,1^2)$, then $\pi$ is graphic.

\par
    {\bf Lemma 2.5 [31]}  Let $\pi=(4^x,3^y,2^z,1^m)$ with even $\sigma(\pi)$, $x+y+z+m=n\geq5$
    and $x\geq1$.
Then $\pi\in GS_n$ if and only if $\pi \not \in A$, where
$A=\{(4,3^2,1^2), (4,3,1^3), (4^2,2,1^2), (4^2,3,2,1), (4^3,1^2),
(4^3,2^2), (4^3,3,1), (4^4,2),$ $(4^2,3,1^3), (4^2,1^4),
(4^3,2,1^2), (4^4,1^2), (4^3,1^4)\}$.\par

    {\bf Lemma 2.6 (Kleitman and Wang [15])}\ \   $\pi$ is
graphic if and only if $\pi_k^\prime$ is graphic.
 \par
    The following corollary is obvious.\par
\par
    {\bf Corollary 2.7}\ \    Let $H$ be a simple graph. If $\pi^\prime$ is
 potentially $H$-graphic, then $\pi$ is
 potentially $H$-graphic.
 \par
 In order to prove our main result, we need the following
 definitions and proposition in Yin [28].
 \par
 Let $s\geq2$ and $\pi=(d_1,d_2,\cdots,d_{s+2},\cdots,d_n)$ be a non-increasing sequence of
    nonnegative integers, where $d_1\leq n-2$ and $d_2\geq s$.
    Denote $$ \rho_s^\prime(\pi)=\left\{
    \begin{array}{ll}(d_2-1,d_3-1,\cdots,d_{s+2}-1,d_{s+3}^{(1)},
    \cdots,d_n^{(1)}), \\ \mbox{ if $d_2\geq {s+1};$}\\
    (d_2,d_3-1,\cdots,d_{s+2}-1,d_{s+3}-1,\cdots,d_{d_1+2}-1,d_{d_1+3},\cdots,d_n),
     \\ \mbox{if $d_2=s,$} \end{array} \right. $$
     where $d_{s+3}^{(1)}\geq\cdots\geq d_n^{(1)}$ is a rearrangement of
     $d_{s+3}-1,\cdots,d_{d_1+1}-1,d_{d_1+2},\cdots,d_n$.\par
     From $\rho_s^\prime(\pi)$, we construct the sequence $$ \rho_s(\pi)=\left\{
    \begin{array}{ll}(d_3-2,\cdots,d_{s+2}-2,d_{s+3}^{(2)},
    \cdots,d_n^{(2)}), \\ \mbox{ if $d_2\geq {s+1};$}\\
    (d_3-2,\cdots,d_{s+2}-2,d_{s+3}-1,\cdots,d_{d_1+2}-1,d_{d_1+3},\cdots,d_n),
     \\ \mbox{if $d_2=s,$} \end{array} \right. $$
     where $d_{s+3}^{(2)}\geq\cdots\geq d_n^{(2)}$ is a rearrangement of
     $d_{s+3}^{(1)}-1,\cdots,d_{d_2+1}^{(1)}-1,d_{d_2+2}^{(1)},\cdots,d_n^{(1)}$.\par

\par
    {\bf Proposition 2.8  (Yin [28])}\ \  Let $s\geq2$ and
    $\pi=(d_1,d_2,\cdots,d_{s+2},\cdots,d_n)$ be a non-increasing sequence of
    nonnegative integers, where $d_1\leq n-2$ and $d_2\geq s$.
    If $\rho_s(\pi)$
    is graphic, then $\pi$ is potentially $K_{2,s}$-graphic.
\par
\section{ Main Theorems} \par
\par
\textbf{\noindent Theorem 3.1}  Let $n\geq7$ and
$\pi=(d_1,d_2,\cdots,d_n)\in GS_n$. Then $\pi$ is potentially
$K_{2,5}$-graphic if and only if the following conditions hold:
\par
  (1) $d_2\geq5$, $d_7\geq2$;\par
  (2) If $d_1=n-1$ and $d_2=5$, then $d_3=5$ and $d_7\geq3$;\par
  (3) $\pi=(n-l,5^i,4^j,3^k,2^t,1^{n-7})$  implies
$(3^{i-1},2^j,1^{k+l-2})$ is graphic, where $n-l\geq5$, $l=2,3,4,
i\geq1$ and
 $i+j+k+t=6$;\par
  (4) $\pi$ is not one of the following sequences:\par
     $(n-1,5^4,3^2,1^{n-7})$,\ \  $(n-1,5^3,3^3,1^{n-7})$,\ \
       $(n-1,5^3,3^4,1^{n-8})$,\ \
      $(n-1,5^2,3^5,1^{n-8})$, \ \
      $(n-1,5^2,3^6,1^{n-9})$,\ \ $(n-2,5^4,2^3,1^{n-8})$,\ \
      $(n-2,5^3,3,2^3,1^{n-8})$,\ \ $(n-2,5^2,2^5,1^{n-8})$, \ \ $(n-2,5^2,4,2^4,1^{n-8})$,\ \
      $(n-3,5^3,2^4,1^{n-8})$.\par
      $n=8:$ \ \  $(6^4,3^4)$, $(6^3,4^2,2^3)$, $(6^3,3^2,2^3)$,
      $(6^2,5,3,2^4)$, $(6^2,4,2^5)$, $(5^3,3,2^4)$,\par
      \ \ \ \ \ \ \ \ \ \ \  $(5^2,4,2^5)$, $(5^2,2^6)$.\par
      $n=9:$ \ \  $(7,6^2,3,2^5)$, $(7,6,5,2^6)$, $(6^3,4,2^5)$, $(6^3,2^6)$,
      $(6,5^2,2^6)$, $(5^2,2^7)$,\par
      \ \ \ \ \ \ \ \ \ \ \
      $(6^3,3,2^4,1)$,  $(6^2,5,2^5,1)$, $(5^3,2^5,1)$.\par
      $n=10:$ \ \ $(8,6^2,2^7)$, $(7^3,3,2^6)$, $(7^2,6,2^7)$,
      $(6^3,2^7)$, $(7,6^2,2^6,1)$, $(6^3,2^5,1^2)$.\par
      $n=11:$ \ \  $(8,7^2,2^8)$, $(7^3,2^7,1)$.\par
      $n=12:$ \ \  $(8^3,2^9)$.
\par
{\bf Proof:} First we show the conditions (1)-(4) are necessary
conditions for $\pi$ to be potentially $K_{2,5}$-graphic. Assume
that $\pi$ is potentially $K_{2,5}$-graphic. $(1)$ is obvious. If
$d_1=n-1$ and $d_2=5$, then the residual sequence
$\pi_1^\prime=(d_2-1,d_3-1,\cdots,d_n-1)$ is potentially
$K_{2,4}$-graphic, and hence $d_2-1=d_3-1=4$ and $d_7-1\geq2$, i.e.,
$d_3=5$ and $d_7\geq3$. Hence, (2) holds. If
$\pi=(n-2,5^i,4^j,3^k,2^t,1^{n-7})$ is potentially
$K_{2,5}$-graphic, then according to theorem 2.1, there exists a
realization $G$ of $\pi$ containing $K_{2,5}$ as a subgraph so that
the vertices of $K_{2,5}$ have the largest degrees of $\pi$.
Therefore, the sequence $\pi_1=(n-7,0,3^{i-1},2^j,1^k,0^t,1^{n-7})$
obtained from $G-K_{2,5}$ is graphic. Since the edges of $K_{2,5}$
have been removed from the realization of $\pi_1$, thus,
$(3^{i-1},2^j,1^k)$ must be graphic. Similarly, with the same
argument as above, one can show that
$\pi=(n-l,5^i,4^j,3^k,2^t,1^{n-7})$  implies
$(3^{i-1},2^j,1^{k+l-2})$ is graphic for the cases $l=3$ and $l=4$.
Hence, (3) holds. Now it is easy to check that $(6^4,3^4)$,
$(6^3,4^2,2^3)$, $(6^3,3^2,2^3)$,
      $(6^2,5,3,2^4)$, $(6^2,4,2^5)$, $(5^3,3,2^4)$,
      $(5^2,4,2^5)$,  $(5^2,2^6)$,
       $(7,6^2,3,2^5)$, $(7,6,5,2^6)$, $(6^3,4,2^5)$, $(6^3,2^6)$,
      $(6,5^2,2^6)$, $(5^2,2^7)$, $(6^3,3,2^4,1)$,  $(6^2,5,2^5,1)$, $(5^3,2^5,1)$, $(8,6^2,2^7)$, $(7^3,3,2^6)$, $(7^2,6,2^7)$,
      $(6^3,2^7)$, $(7,6^2,2^6,1)$, $(6^3,2^5,1^2)$, $(8,7^2,2^8)$, $(7^3,2^7,1)$ and $(8^3,2^9)$
       are not potentially $K_{2,5}$-graphic. Since
       $\pi_1^\prime=(4^4,2^2)$, $(4^3,2^3)$, $(4^3,2^4)$, $(4^2,2^5)$
and $(4^2,2^6)$ are not potentially $K_{2,4}$(by theorem 2.3) or
$K_{1,5}$-graphic, we have $\pi\neq(n-1,5^4,3^2,1^{n-7})$,
$(n-1,5^3,3^3,1^{n-7})$, $(n-1,5^3,3^4,1^{n-8})$,
$(n-1,5^2,3^5,1^{n-8})$ and $(n-1,5^2,3^6,1^{n-9})$. If
$\pi=(n-2,5^4,2^3,1^{n-8})$ is potentially $K_{2,5}$-graphic, then
according to theorem 2.1, there exists a realization $G$ of $\pi$
containing $K_{2,5}$ as a subgraph so that the vertices of $K_{2,5}$
have the largest degrees of $\pi$. Therefore, the sequence
$\pi_1=(n-7,0,3^3,0^2,2,1^{n-8})$ obtained from $G-K_{2,5}$ must be
graphic. Since the edges of $K_{2,5}$ have been removed from the
realization of $\pi_1$,  $\pi_2=(3^3,1)$ is graphic, a
contradiction. Thus, $\pi\neq(n-2,5^4,2^3,1^{n-8})$. Similarly, one
can show that $\pi\neq(n-2,5^3,3,2^3,1^{n-8})$,
$(n-2,5^2,2^5,1^{n-8})$, $(n-2,5^2,4,2^4,1^{n-8})$ and
$(n-3,5^3,2^4,1^{n-8})$. Hence, (4) holds.
\par
 Next, we will prove the sufficient conditions. Suppose $\pi=(d_1,d_2,\cdots,d_n)\in GS_n$
 satisfies the conditions (1)-(4).\par
 If $d_1=n-1$, consider the residual sequence
 $\pi_1^\prime=(d_1^\prime,d_2^\prime,\cdots,d_{n-1}^\prime)$
 obtained by laying off $d_1$ form $\pi$. If $d_2\geq6$, then
 $d_1^\prime=d_2-1\geq5$. Thus, $\pi_1^\prime$ is potentially
 $K_{1,5}$-graphic and so $\pi$ is
potentially $K_{2,5}$-graphic. If $d_2=5$, by $\pi$ satisfies
condition (2), we have $d_1^\prime=d_2-1=4$, $d_2^\prime=d_3-1=4$
and $d_6^\prime=d_7-1\geq2$. Since $\pi$ satisfies condition (4),
$\pi_1^\prime\neq(4^4,2^2)$, $(4^3,2^3)$, $(4^3,2^4)$, $(4^2,2^5)$
and $(4^2,2^6)$. If $\pi_1^\prime \not \in A$, where
$A=\{(5^2,4,2^4),$ $(5^3,3,2^3)$, $(6,5^2,2^5)$, $(5^3,2^4,1)$,
$(6^3,2^6)$, $(n-2,4^2,3^4,1^{n-8})$, $(n-2,4^2,3^5,1^{n-9})$,
$(n-3,4^2,2^3,1^{n-7})(n\geq8)$, $(n-3,4^3,2^2,1^{n-7})(n\geq8)$\},
then $\pi_1^\prime$ is potentially
 $K_{2,4}$-graphic by theorem 2.3 and so $\pi$ is
potentially $K_{2,5}$-graphic. If $\pi_1^\prime  \in A$, then
$\pi_1^\prime$ is potentially
 $K_{1,5}$-graphic, thus $\pi$ is
potentially $K_{2,5}$-graphic. Suppose $d_1\leq n-2$. \par
 Our proof is by induction on $n$. We first prove the base case
where $n=7$. In this case, $d_1=d_2=5$, i.e.,
$\pi=(5^i,4^j,3^k,2^{7-i-j-k})$ where $i\geq2$. Then
$\rho_5(\pi)=(3^{i-2},2^j,1^k,0^{7-i-j-k})$. By $\pi$ satisfies (3),
$\rho_5(\pi)$ is graphic, and so $\pi$ is potentially
$K_{2,5}$-graphic by proposition 2.8. Now suppose that the
sufficiency holds for $n-1(n\geq8)$, we will show that $\pi$ is
potentially $K_{2,5}$-graphic in terms of the following cases:
\par
\textbf{Case 1:} $d_n\geq5$.  Consider
$\pi^\prime=(d_1^\prime,d_2^\prime,\cdots,d_{n-1}^\prime)$. Clearly,
$\pi^\prime$ satisfies (1) and (4). If $\pi^\prime$ also satisfies
(2)-(3), then by the induction hypothesis, $\pi^\prime$ is
potentially $K_{2,5}$-graphic, and hence so is $\pi$.
\par
  If $d_1^\prime=n-2$ and $d_2^\prime=5$, by $d_1\leq n-2$, then
  $d_1=d_2=\cdots=d_6=n-2$ and $n=8$. Thus, $\pi^\prime=(6,5^6)$
  which satisfies condition (2).
\par
If $\pi^\prime=(n-1-l,5^i,4^{6-i})$, then $l=2$ and $n=8$. Thus,
$\pi^\prime=(5^2,4^5)$, $(5^4,4^3)$ or $(5^6,4)$. Since $(2^5)$,
$(3^2,2^3)$ and $(3^4,2)$ are graphic, $\pi^\prime$ satisfies
condition (3).
\par
\textbf{Case 2:} $d_n=4$. Consider
$\pi^\prime=(d_1^\prime,d_2^\prime,\cdots,d_{n-1}^\prime)$ where
$d_{n-3}^\prime\geq4$ and $d_{n-1}^\prime\geq3$. Clearly,
$\pi^\prime$ satisfies (4). If $\pi^\prime$ also satisfies (1)-(3),
then by the induction hypothesis, $\pi^\prime$ is potentially
$K_{2,5}$-graphic, and hence so is $\pi$.
 \par
If $\pi^\prime$ does not satisfy $(1)$, i.e., $d_2^\prime=4$, then
$d_2=5$. We will proceed with the following two cases: $d_1=5$ and
$d_1\geq6$.
\par
\textbf{Subcase 1:} $d_1=5$. Then $\pi=(5^k,4^{n-k})$ where $2\leq
k\leq5$. Since $\sigma(\pi)$ is even, we have $k=2$ or $k=4$, i.e.,
$\pi=(5^2,4^{n-2})$ or $(5^4,4^{n-4})$. Then
$\rho_5(\pi)=(2^5,4^{n-7})$ or $(3^2,2^3,4^{n-7})$. By lemma 2.5,
$\rho_5(\pi)$ is graphic, and so $\pi=(5^k,4^{n-k})$ is potentially
$K_{2,5}$-graphic by proposition 2.8.\par \textbf{Subcase 2:}
$d_1\geq6$. Then $\pi=(d_1,5^k,4^{n-1-k})$ where $1\leq k\leq3$,
and, $d_1$ and $k$ have the same parity.
\par
If $k=1$, then $\pi=(d_1,5,4^{n-2})$ and
$\rho_5(\pi)=(2^5,3^{d_1-5},4^{n-2-d_1})$. By lemma 2.4 and lemma
2.5, $\rho_5(\pi)$ is graphic, and so $\pi=(d_1,5,4^{n-2})$ is
potentially $K_{2,5}$-graphic by proposition 2.8.
\par
If $k=2$, then $\pi=(d_1,5^2,4^{n-3})$ and
$\rho_5(\pi)=(3,2^4,3^{d_1-5},4^{n-2-d_1})$. By lemma 2.4 and lemma
2.5, $\rho_5(\pi)$ is graphic, and so $\pi=(d_1,5^2,4^{n-3})$ is
potentially $K_{2,5}$-graphic by proposition 2.8. Similarly, with
the same argument as above, one can show that
$\pi=(d_1,5^3,4^{n-4})$ is also potentially $K_{2,5}$-graphic.
\par
If $d_1^\prime=n-2$ and $d_2^\prime=5$, by $d_1\leq n-2$, then
$d_1=d_2=d_3=d_4=d_5=n-2$ and $n=8$. Thus, $\pi^\prime=(6,5^6)$ or
$(6,5^4,4^2)$. Therefore, $\pi^\prime$ satisfies condition (2).
\par
If $\pi^\prime=(n-1-l,5^i,4^j,3^{6-i-j})$, then $n=8$, $l=2$ and
$i+j\geq4$. Hence, $\pi^\prime$ is one of the following: $(5^6,4)$,
$(5^4,4^3)$, $(5^2,4^5)$, $(5^4,4,3^2)$, $(5^2,4^3,3^2)$,
$(5^5,4,3)$, $(5^3,4^3,3)$. By lemma 2.4, $(3^4,2)$, $(3^2,2^3)$,
$(2^5)$, $(3^2,2,1^2)$, $(2^3,1^2)$, $(3^3,2,1)$ and $(3,2^3,1)$ are
graphic. Thus, $\pi^\prime$ satisfies condition (3).
\par
 \textbf{Case 3:} $d_n=3$.
Consider $\pi^\prime=(d_1^\prime,d_2^\prime,\cdots,d_{n-1}^\prime)$
where $d_2^\prime\geq4$, $d_{n-2}^\prime\geq3$ and
$d_{n-1}^\prime\geq2$. If $\pi^\prime$ satisfies (1)-(4), then by
the induction hypothesis, $\pi^\prime$ is potentially
$K_{2,5}$-graphic, and hence so is $\pi$.
\par
If $\pi^\prime$ does not satisfy $(1)$, i.e., $d_2^\prime=4$, then
$d_2=5$. We will proceed with the following two cases: $d_1=5$ and
$d_1\geq6$.
\par
\textbf{Subcase 1:} $d_1=5$. Then $\pi=(5^i,4^j,3^{n-i-j})$ where
$2\leq i\leq4$, $n-i-j\geq1$ and $n-j$ is even. If $i=2$, i.e.,
$\pi=(5^2,4^j,3^{n-2-j})$. Then
$\rho_5(\pi)=(2^5,4^{j-5},3^{n-2-j})$$(j\geq5)$ or
$\rho_5(\pi)=(2^j,1^{5-j},3^{n-7})$$(j<5)$. By lemma 2.4 and lemma
2.5, $\rho_5(\pi)$ is graphic, and so $\pi=(5^2,4^j,3^{n-2-j})$ is
potentially $K_{2,5}$-graphic by proposition 2.8. Similarly, with
the same argument as above, one can show that
$\pi=(5^i,4^j,3^{n-i-j})$ is also potentially $K_{2,5}$-graphic for
the cases $i=3$ and $i=4$.
\par
\textbf{Subcase 2:} $d_1\geq6$. Then $\pi=(d_1,5^i,4^j,3^{n-1-i-j})$
where $1\leq i\leq2$, $n-1-i-j\geq1$, and, $d_1$ and $n-1-j$ have
the same parity. If $i=1$, then $\pi=(d_1,5,4^j,3^{n-2-j})$. If
$j<5$, then $\rho_5(\pi)=(2^j,1^{5-j},2^{d_1-5},3^{n-2-d_1})$. If
$5\leq j<d_1$, then
$\rho_5(\pi)=(2^5,3^{j-5},2^{d_1-j},3^{n-2-d_1})$. If $j\geq d_1$,
then $\rho_5(\pi)=(2^5,3^{d_1-5},4^{j-d_1},3^{n-2-j})$. By lemma 2.4
and lemma 2.5, $\rho_5(\pi)$ is graphic, and so
$\pi=(d_1,5^i,4^j,3^{n-1-i-j})$ is potentially $K_{2,5}$-graphic by
proposition 2.8. Similarly, with the same argument as above, one can
show that $\pi=(d_1,5^2,4^j,3^{n-3-j})$ is also potentially
$K_{2,5}$-graphic.
\par
If $d_1^\prime=n-2$ and $d_2^\prime=5$, by $d_1\leq n-2$, then
$d_1=d_2=d_3=d_4=n-2$ and $n=8$. In this case,
$\pi^\prime=(6,5^3,d_5^\prime,d_6^\prime,d_7^\prime)$ where $3\leq
d_7^\prime\leq d_6^\prime\leq d_5^\prime\leq 5$ and
$\sigma(\pi^\prime)$ is even. Clearly, $\pi^\prime$ satisfies
condition (2).
\par
  If $\pi^\prime=(n-1-l,5^i,4^j,3^k,2^{6-i-j-k})$, then $n=8,$ $l=2$ and
  $i+j+k\geq5$. If $\pi^\prime=(5,5^i,4^j,3^{6-i-j})$, then it is
  easy to see that $(3^{i-1},2^j,1^{6-i-j})$ is graphic by lemma
  2.4(since $(i-1)+j+(6-i-j)=5>4$). If $\pi^\prime=(5,5^i,4^j,3^{5-i-j},2)$, then
  $d_3=\cdots=d_7=3$. We have $\pi^\prime=(5^2,3^4,2)$. It follows
  $(3^{i-1},2^j,1^{5-i-j})=(1^4)$, which is also graphic. In other
  words, $\pi^\prime$ satisfies (3).
\par
  If $\pi^\prime$ does not satisfy $(4)$, since $d_1\leq n-2$ and $\pi\neq(6^4,3^4)$,
 then $\pi^\prime=(6,5^4,3^2)$ or $(6^4,3^4)$. Thus,
  $\pi=(6^4,5,3^3)$ or $(7^3,6,3^5)$. It is easy to check
  that both of them are potentially
$K_{2,5}$-graphic.
\par
 \textbf{Case 4:} $d_n=2$.
Consider $\pi^\prime=(d_1^\prime,d_2^\prime,\cdots,d_{n-1}^\prime)$
where $d_2^\prime\geq4$ and $d_{n-1}^\prime\geq2$. If $\pi^\prime$
satisfies (1)-(4), then by the induction hypothesis, $\pi^\prime$ is
potentially $K_{2,5}$-graphic, and hence so is $\pi$.
\par
If $\pi^\prime$ does not satisfy $(1)$, i.e., $d_2^\prime=4$, then
$d_2=5$. We will proceed with the following two cases: $d_1=5$ and
$d_1\geq6$.
\par
\textbf{Subcase 1:} $d_1=5$. Then $\pi=(5^k,4^i,3^j,2^{n-k-i-j})$
where $2\leq k\leq3$, $n-k-i-j\geq1$ and $k+j$ is even. If $k=2$,
i.e., $\pi=(5^2,4^i,3^j,2^{n-2-i-j})$, then
$\rho_5(\pi)=(d_3-2,\cdots,d_7-2,d_8,\cdots,d_n)$. If
$m(\rho_5(\pi))=4$, then $i\geq 6$ and so
$\rho_5(\pi)=(2^5,4^{i-5},3^j,2^{n-2-i-j})$. By lemma 2.5,
$\rho_5(\pi)$ is graphic. If $m(\rho_5(\pi))=3$, then $i\leq5$,
$i+j\geq6$ and so $\rho_5(\pi)=(2^i,1^{5-i},3^{i+j-5},2^{n-2-i-j})$,
it follows from lemma 2.4 that $\rho_5(\pi)$ is also graphic. If
$m(\rho_5(\pi))=2$, then $i+j\leq5$ and so
$\rho_5(\pi)=(2^i,1^j,0^{5-i-j},2^{n-7})$. In this case,
$\rho_5(\pi)$ is not graphic if and only if $\rho_5(\pi)=(2)$ or
$(2^2)$ which is impossible since $\pi\neq(5^2,2^6)$, $(5^2,2^7)$
and $(5^2,4,2^5)$. Thus, $\pi=(5^2,4^i,3^j,2^{n-2-i-j})$ is
potentially $K_{2,5}$-graphic by proposition 2.8. Similarly, one can
show that $\pi=(5^3,4^i,3^j,2^{n-3-i-j})$ is also potentially
$K_{2,5}$-graphic.
\par
\textbf{Subcase 2:} $d_1\geq6$. Then
$\pi=(d_1,5,4^i,3^j,2^{n-2-i-j})$ where $n-2-i-j\geq1$, and, $d_1$
and $j+1$ have the same parity. Thus,
$\rho_5(\pi)=(d_3-2,\cdots,d_7-2,d_8-1,\cdots,d_{d_1+2}-1,d_{d_1+3},\cdots,
d_n)$. If $i+j<5$, then
$\rho_5(\pi)=(2^i,1^j,0^{5-i-j},1^{d_1-5},2^{n-2-d_1})$.  If $i\geq
d_1$, then $\rho_5(\pi)=(2^5,3^{d_1-5},4^{i-d_1},3^j,2^{n-2-i-j})$.
If $5\leq i<d_1$ and $i+j\geq d_1$, then
$\rho_5(\pi)=(2^5,3^{i-5},2^{d_1-i},3^{i+j-d_1},2^{n-2-i-j})$. If
$5\leq i+j<d_1$, then $\rho_5(\pi)=(2^5,3^{i-5},$
 $2^j,1^{d_1-i-j},2^{n-2-d_1})$$(i\geq5)$ or
$\rho_5(\pi)=(2^i,1^{5-i},2^{i+j-5},1^{d_1-i-j},2^{n-2-d_1})$
$(i<5)$. By lemma  2.4 and lemma 2.5, in the above cases,
$\rho_5(\pi)$ is graphic, and so $\pi$ is potentially
$K_{2,5}$-graphic by proposition  2.8.
\par
If $d_1^\prime=n-2$ and $d_2^\prime=5$, by $d_1\leq n-2$, we have
$d_1=d_2=d_3=n-2$ and $n=8$. If $d_7\geq3$, then $\pi^\prime$
satisfies (2). If $d_7=2$, then $\pi=(6^3,d_4,d_5,d_6,2^2)$ where
$2\leq d_6\leq d_5\leq d_4\leq5$. Since $\pi\neq(6^3,4^2,2^3)$ and
$(6^3,3^2,2^3)$, then $\pi=(6^3,5^2,4,2^2)$, $(6^3,5,4,3,2^2)$,
$(6^3,4^3,2^2)$ or $(6^3,4,3^2,2^2)$. It is easy to check that all
of these are potentially $K_{2,5}$-graphic.
\par
If $\pi^\prime=(n-1-l,5^i,4^j,3^k,2^{6-i-j-k})$, then $n=8$ and
$l=2$, i.e., $\pi^\prime=(5,5^i,4^j,3^k,2^{6-i-j-k})$. If
$(3^{i-1},2^j,1^k)$ is graphic, then $\pi^\prime$ satisfies (3). If
$(3^{i-1},2^j,1^k)$ is not graphic, then $(3^{i-1},2^j,1^k)$ $\in$
$\{(3^3,1),(3^2,1^2),(3^2,2),$  $(3,2,1),(3,1),(2^2),(2)\}$. By
$\pi\neq(6^2,5,3,2^4)$, $(6^2,4,2^5)$, $(6,5^2,2^5)$, then
$\pi^\prime$ is one of the following: $(5^5,3,2)$, $(5^4,3^2,2)$,
$(5^4,4,2^2)$, $(5^3,4,3,2^2)$, $(5^2,4^2,2^3)$. Since
$\pi\neq(6,5^4,2^3)$,$(6,5^3,3,2^3)$, $(6,5^2,4,2^4)$ and
$(5^4,2^4)$, then $\pi$ is one of the following: $(6^2,5^3,3,2^2)$,
$(6^2,5^2,3^2,2^2)$, $(6^2,5^2,4,2^3)$, $(6^2,5,4,3,2^3)$,
$(6^2,4^2,2^4)$. It is easy to check that all of these are
potentially $K_{2,5}$-graphic.
\par
If $\pi^\prime$ does not satisfy (4), since $d_1\leq n-2$ and
$\pi\neq(7,6,5,2^6)$, $(6^3,2^6)$, $(7^2,6,2^7)$, then $\pi^\prime$
is one of the following: \par $n-1=7:$  $(6,5^4,3^2)$,
$(6,5^3,3^3)$,
 \par
 $n-1=8:$ $(6,5^4,2^3)$,
$(6,5^3,3,2^3)$, $(6,5^2,4,2^4)$, $(5^4,2^4)$, $(6^4,3^4)$,
$(6^3,4^2,2^3)$, \par  \ \ \ \ \ \ \ \ \ \ \ \ $(6^3,3^2,2^3)$,
      $(6^2,5,3,2^4)$, $(6^2,4,2^5)$, $(5^3,3,2^4)$,
      $(5^2,4,2^5)$,  $(5^2,2^6)$,\par
      $n-1=9:$  $(7,6^2,3,2^5)$, $(7,6,5,2^6)$, $(6^3,4,2^5)$,
      $(6,5^2,2^6)$, $(5^2,2^7)$,  \par
            $n-1=10:$ $(8,6^2,2^7)$, $(7^3,3,2^6)$, $(7^2,6,2^7)$,
      $(6^3,2^7)$,\par
      $n-1=11:$  $(8,7^2,2^8)$,\par
      $n-1=12:$  $(8^3,2^9)$,\par
      Since $\pi\neq(6^3,4,2^5)$, $(7,6^2,3,2^5)$, $(6,5^2,2^6)$,
      $(7^3,3,2^6)$, $(8,6^2,2^7)$, $(7^2,6,2^7)$, $(6^3,2^7)$,
      $(8,7^2,2^8)$, $(8^3,2^9)$, then $\pi$ is one of the
following:
\par
 $n=8:$ $(6^3,5^2,3^2,2)$,
$(6^3,5,3^3,2)$.
\par
 $n=9:$ $(7,6,5^3,2^4)$, $(6^3,5^2,2^4)$, $(7,6,5^2,3,2^4)$, $(6^3,5,3,2^4)$,
 $(7,6,5,4,2^5)$, \par
 \ \ \ \ \ \ \ \ \ \ $(6^2,5^2,2^5)$,  $(7^2,6^2,3^4,2)$,
$(7^2,6,4^2,2^4)$, $(7^2,6,3^2,2^4)$, $(7^2,5,3,2^5)$, \par
 \ \ \ \ \ \ \ \ \ \
$(7^2,4,2^6)$, $(6^2,5,3,2^5)$, $(6^2,4,2^6)$, $(6^2,2^7)$.
\par
 $n=10:$ $(8,7,6,3,2^6)$, $(8,7,5,2^7)$,
$(7^2,6,4,2^6)$, $(7,6,5,2^7)$, $(6^2,2^8)$.
\par
 $n=11:$ $(9,7,6,2^8)$, $(8^2,7,3,2^7)$, $(8^2,6,2^8)$,
 $(7^2,6,2^8)$.

\par
 $n=12:$ $(9,8,7,2^9)$.
 \par
 $n=13:$ $(9^2,8,2^{10})$.\par
  It is easy to check that all of
these are potentially $K_{2,5}$-graphic.
\par
 \textbf{Case 5:} $d_n=1$.
Consider $\pi^\prime=(d_1^\prime,d_2^\prime,\cdots,d_{n-1}^\prime)$
where $d_1^\prime\geq5$, $d_2^\prime\geq4$ and $d_7^\prime\geq2$. If
$\pi^\prime$ satisfies (1)-(4), then by the induction hypothesis,
$\pi^\prime$ is potentially $K_{2,5}$-graphic, and hence so is
$\pi$.
\par
If $\pi^\prime$ does not satisfy $(1)$, i.e., $d_2^\prime=4$. Then
$d_1=d_2=5$ and $d_3\leq4$, i.e., $\pi=(5^2,d_3,\cdots,d_n)$ where
$d_7\geq2$ and $d_n=1$. Consider
$\rho_5(\pi)=(d_3-2,\cdots,d_7-2,d_8,\cdots,d_n)$. If
$m(\rho_5(\pi))=4$, then $d_3=\cdots=d_8=4$, i.e.,
$\rho_5(\pi)=(2^5,4,d_9,\cdots,d_{n-1},1)$. By lemma 2.5,
$\rho_5(\pi)$ is graphic. If $m(\rho_5(\pi))=3$, then $d_8=3$, i.e.,
$\rho_5(\pi)=(d_3-2,\cdots,d_7-2,3,d_9,\cdots,d_{n-1},1)$ where
$1\leq d_7-2\leq d_6-2\leq \cdots \leq d_3-2\leq2$. By lemma 2.4,
$\rho_5(\pi)$ is also graphic. If $1\leq m(\rho_5(\pi))\leq2$, it
follows from $h(\rho_5(\pi))=1$ and theorem 2.2  that $\rho_5(\pi)$
is graphic. Hence, $\pi$ is potentially $K_{2,5}$-graphic by
proposition 2.8.
\par
If $d_1^\prime=n-2$ and $d_2^\prime=5$, by $d_1\leq n-2$, we have
$d_1=d_2=n-2$ and $n=8$, i.e., $\pi=(6^2,d_3,\cdots,d_7,1)$. If
$d_3=5$ and $d_7\geq3$, then $\pi^\prime$ satisfies (2). If
$d_3\leq4$, then $\pi$ is one of the following: $(6^2,4^4,3,1)$,
$(6^2,4^3,3,2,1)$, $(6^2,4^2,3^3,1)$, $(6^2,4^2,3,2^2,1)$,
$(6^2,4,3^3,2,1)$, $(6^2,4,3,2^3,1)$, $(6^2,3^5,1)$,
$(6^2,3^3,2^2,1)$, $(6^2,3,2^4,1)$. If $d_3=5$ and $d_7=2$, then
$\pi$ is one of the following: $(6^2,5^3,4,2,1)$,
$(6^2,5^2,4,3,2,1)$, $(6^2,5,4^3,2,1)$, $(6^2,5,4^2,2^2,1)$,
$(6^2,5,4,3^2,2,1)$, $(6^2,5,3^2,2^2,1)$. It is easy to check that
all of the above sequences are potentially $K_{2,5}$-graphic.
\par
If $\pi^\prime=(n-1-l,5^i,4^j,3^k,2^t,1^{n-8})$, then there are
three subcases:
\par
\textbf{Subcase 1: } $n-1-l=5$. If $j=0$, then
$\pi^\prime=(5,5^i,3^k,2^t,1^{n-8})$ and
$\pi=(6,5^i,3^k,2^t,1^{n-7})$. By $\pi$ satisfies (3), then
$(3^{i-1},1^{k+n-8})$ is graphic. Thus, $\pi^\prime$ satisfies (3).
If $j\geq1$, then $\pi=(6,5^i,4^j,3^k,2^t,1^{n-7})$ or
$(5^{i+2},4^{j-1},3^k,2^t,1^{n-7})$. If
$\pi=(6,5^i,4^j,3^k,2^t,1^{n-7})$, then  with the same argument as
above, one can show that $\pi^\prime$ satisfies (3). If
$\pi=(5^{i+2},4^{j-1},3^k,2^t,1^{n-7})$, then
$\rho_5(\pi)=(3^i,2^{j-1},1^k,0^t,1^{n-7})$. Since $\pi$ satisfies
(3), $\rho_5(\pi)$ is graphic. Thus, $\pi$ is potentially
$K_{2,5}$-graphic by proposition 2.8.
\par
\textbf{Subcase 2: } $n-1-l=6$. Then
$\pi=(7,5^i,4^j,3^k,2^t,1^{n-7})$ or
$(6^2,5^{i-1},4^j,3^k,2^t,1^{n-7})$. If
$\pi=(7,5^i,4^j,3^k,2^t,1^{n-7})$, with the same argument as
subcase1, we have $\pi^\prime$ satisfies (3). If
$\pi=(6^2,5^{i-1},4^j,3^k,2^t,1^{n-7})$, then
$\rho_5(\pi)=(3^{i-1},2^j,1^k,0^t,1^{n-7})$ where $n\geq9$. In this
case, $\rho_5(\pi)$ is not graphic if and only if
$\rho_5(\pi)=(3^2,1^2)$ which is impossible since
$\pi=(6^2,5^2,2^3,1^2)$ is not graphic. Thus,
$\pi=(6^2,5^{i-1},4^j,3^k,2^t,1^{n-7})$ is potentially
$K_{2,5}$-graphic by proposition 2.8.
\par
\textbf{Subcase 3: } $n-1-l\geq 7$. Then
$\pi=(n-l,5^i,4^j,3^k,2^t,1^{n-7})$. By $\pi$ satisfies (3), then
$(3^{i-1},2^j,1^{k+l-2})$ is graphic. In other words, $\pi^\prime$
satisfies (3).
\par
If $\pi^\prime$ does not satisfy (4), since
$\pi\neq(n-1,5^3,3^4,1^{n-8})$, $(n-1,5^2,3^5,1^{n-8})$,
$(n-1,5^2,3^6,1^{n-9})$, $(n-2,5^2,2^5,1^{n-8})$,
$(n-3,5^3,2^4,1^{n-8})$, $(n-3,5^2,3,2^4,$   $1^{n-8})$,
$(6^2,5,2^5,1)$, $(7,6^2,2^6,1)$,  then $\pi^\prime$ is one of the
following: \par $n-1=7:$ $(6,5^4,3^2)$, $(6,5^3,3^3)$,
\par $n-1=8:$ $(6,5^4,2^3)$, $(6,5^3,3,2^3)$, $(6,5^2,4,2^4)$, $(6^4,3^4)$,
$(6^3,4^2,2^3)$, \par \ \ \ \ \ \ \ \ \ \ \ $(6^3,3^2,2^3)$,
      $(6^2,5,3,2^4)$, $(6^2,4,2^5)$, $(5^3,3,2^4)$, $(5^2,4,2^5)$,  $(5^2,2^6)$,
\par $n-1=9:$ $(7,6^2,3,2^5)$, $(7,6,5,2^6)$, $(6^3,4,2^5)$,
      $(6,5^2,2^6)$, $(5^2,2^7)$, \par
      \ \ \ \ \ \ \ \ \ \ \ \ \ \ \ $(6^3,3,2^4,1)$,  $(6^2,5,2^5,1)$, $(5^3,2^5,1)$,
\par $n-1=10:$ $(8,6^2,2^7)$, $(7^3,3,2^6)$, $(7^2,6,2^7)$,
      $(6^3,2^7)$, $(7,6^2,2^6,1)$,\par
      \ \ \ \ \ \ \ \ \ \ \ \ \ \ \ \  $(6^3,2^5,1^2)$,
\par $n-1=11:$ $(8,7^2,2^8)$, $(7^3,2^7,1)$,
\par $n-1=12:$  $(8^3,2^9)$.\par Since
$\pi\neq(n-1,5^4,3^2,1^{n-7})$, $(n-1,5^3,3^3,1^{n-7})$,
$(n-2,5^4,2^3,1^{n-8})$, $(n-2,5^3,3,2^3,1^{n-8})$,
$(n-2,5^2,4,2^4,1^{n-8})$, $(6^3,3,2^4,1)$, $(5^3,2^5,1)$,
$(6^3,2^5,1^2)$ and $(7^3,2^7,1)$, then $\pi$ is one of the
following:
\par
$n=8:$ $(6^2,5^3,3^2,1)$, $(6^2,5^2,3^3,1)$,
\par
$n=9:$ $(6^2,5^3,2^3,1)$, $(6^2,5^2,3,2^3,1)$, $(6^2,5,4,2^4,1)$,
$(7,6^3,3^4,1)$, \par \ \ \ \ \ \ \ \ \ \ $(7,6^2,4^2,2^3,1)$,
$(7,6^2,3^2,2^3,1)$, $(7,6,5,3,2^4,1)$, $(7,6,4,2^5,1)$,  \par \ \ \
\ \ \ \ \ \ \ $(6,5^2,3,2^4,1)$, $(6,5,4,2^5,1)$, $(6,5,2^6,1)$,
\par
$n=10:$ $(8,6^2,3,2^5,1)$, $(7^2,6,3,2^5,1)$, $(8,6,5,2^6,1)$,
$(7^2,5,2^6,1)$,
 \par \ \ \ \ \ \ \ \ \ \ \
 $(7,6^2,4,2^5,1)$, $(7,5^2,2^6,1)$, $(6^2,5,2^6,1)$, $(6,5,2^7,1)$,
 \par \ \ \ \ \ \ \ \ \ \ \ $(7,6^2,3,2^4,1^2)$, $(7,6,5,2^5,1^2)$,
 $(6,5^2,2^5,1^2)$,
\par
$n=11:$ $(9,6^2,2^7,1)$, $(8,7^2,3,2^6,1)$,  $(8,7,6,2^7,1)$,
$(7,6^2,2^7,1)$, \par \ \ \ \ \ \ \ \ \ \ \ \ $(8,6^2,2^6,1^2)$,
 $(7^2,6,2^6,1^2)$, $(7,6^2,2^5,1^3)$,
\par
$n=12:$ $(9,7^2,2^8,1)$,
 $(8^2,7,2^8,1)$, $(8,7^2,2^7,1^2)$,
\par
$n=13:$ $(9,8^2,2^9,1)$. \par
It is easy to check that all of the
above sequences are potentially $K_{2,5}$-graphic.
\par
\par
\vspace{0.5cm}

\section{  Application }

\par
In the remaining of this section, we will use theorem 3.1 to find
the exact value of $\sigma(K_{2,5},n)$. Note that the value of
$\sigma(K_{2,5},n)$ was determined by Yin et al. in $[26]$ so
another proof is given here.
\par
  \textbf{Theorem }  (Yin et al. [26])  If $n\geq37$, then
    $$ \sigma(K_{2,5},n)=\left\{
    \begin{array}{ll} 5n-3, \ \mbox{ if $n$ is odd,}\\
    5n-2,
     \ \  \mbox{if $n$ is even.} \end{array} \right. $$
\par
\textbf{Proof:} First we claim that for $n\geq37$,
$$ \sigma(K_{2,5},n)\geq \left\{
    \begin{array}{ll} 5n-3, \ \mbox{ if $n$ is odd,}\\
    5n-2,
     \ \  \mbox{if $n$ is even.} \end{array} \right. $$ If $n$ is odd, take $\pi_1=((n-1),5,4^{n-3},3)$, then
$\sigma(\pi_1)=5n-5$, and it is easy to see that $\pi_1$ is not
potentially $K_{2,5}$-graphic by theorem 3.1. If $n$ is even, take
$\pi_2=(n-1,5,4^{n-2})$, then $\sigma(\pi_2)=5n-4$ and  $\pi_2$ is
not potentially $K_{2,5}$-graphic by theorem 3.1. Thus, $$
\sigma(K_{2,5},n)\geq \left\{
    \begin{array}{ll} \sigma(\pi_1)+2=5n-3, \ \mbox{ if $n$ is odd,}\\
    \sigma(\pi_2)+2=5n-2,
     \ \  \mbox{if $n$ is even.} \end{array} \right. $$
\par
  Now we show that if $\pi$ is an $n$-term $(n\geq37)$ graphical
sequence with $\sigma(\pi)\geq5n-3$, then there exists a realization
of $\pi$ containing $K_{2,5}$. Hence, it suffices to show that $\pi$
is potentially $K_{2,5}$-graphic.
\par
If $d_2\leq4$, then $\sigma(\pi)\leq d_1+4(n-1)\leq
n-1+4(n-1)=5n-5<5n-3$, a contradiction. Hence, $d_2\geq5$.
\par
 If $d_7=1$, then $\sigma(\pi)=
d_1+d_2+d_3+d_4+d_5+d_6+(n-6)\leq30+(n-6)+(n-6)=2n+18<5n-3$, a
contradiction. Hence, $d_7\geq2$.
\par
If $d_1=n-1$, $d_2=5$ and $d_3\leq4$, then $\sigma(\pi)\leq
(n-1)+5+4(n-2)=5n-4<5n-3$, a contradiction. If $d_1=n-1$, $d_2=5$
and $d_7\leq2$, then $\sigma(\pi)\leq
(n-1)+5\times5+2(n-6)=3n+12<5n-3$, a contradiction. Hence, $\pi$
satisfies condition (2) in theorem 3.1.
\par
  Since \ \  $\sigma(\pi)\geq5n-3$, it is easy to check that $\pi$ satisfies condition (4) in theorem 3.1.
  Therefore, $\pi$ is potentially $K_{2,5}$-graphic.

\par

\end{document}